\documentclass[12pt]{amsart}

\usepackage[utf8]{inputenc}
\usepackage[english]{babel}

\usepackage[letterpaper,top=2cm,bottom=2cm,left=3cm,right=3cm,marginparwidth=1.75cm]{geometry}

\usepackage{amsmath}
\usepackage{graphicx}
\usepackage[colorlinks=true, allcolors=blue]{hyperref}

\usepackage{amssymb} 
\usepackage{comment}
\usepackage{indentfirst}
\usepackage{csquotes}
\usepackage{theoremref}
\usepackage{hyperref}
\usepackage{geometry}

\hypersetup{hidelinks, colorlinks=true, allcolors=cyan, pdfstartview=Fit, breaklinks=true}

\usepackage{amsthm}

\theoremstyle{plain}

\newtheorem*{corollary*}{Corollary}
\newtheorem*{Example*}{Example}

\newtheorem{lemma}{Lemma}

\theoremstyle{definition}

\numberwithin{definition}{section}
\newtheorem*{def*}{Definition}
\newtheorem*{theorem*}{Theorem}

\newtheorem*{definition*}{Definition}

\theoremstyle{remark}

\newcommand{\bracket}[1]{\left(#1\right)}
\newcommand{\largebracket}[1]{\left\{#1\right\}}

\newcommand{\ceil}[1]{\left\lceil#1\right\rceil}
\newcommand{\modulo}[3]{#1\equiv#2\ \bracket{\mathrm{mod}\ #3}}
\newcommand{\notmodulo}[3]{#1\not\equiv #2\ \bracket{\mathrm{mod}\ #3}}
\newcommand{\emptyfrac}[2]{\genfrac{}{}{0pt}{}{#1}{#2}}
\newcommand{\bbZ}[0]{\mathbb Z}

\title{\textbf{Repeated values of some restricted divisor functions}}
\author{QI-YANG ZHENG}
\date{} 

\address{Department of Mathematics, Sun Yat-sen University(Zhuhai Campus), Zhuhai}
\email{zhengqy29@mail2.sysu.edu.cn}

\begin{document}
\maketitle

\begin{abstract}
We prove that $d_k(n)=d_k(n+B)$ infinitely often for any positive integers $k$ and $B$, where $d_k(n)$ denotes the number of divisors of $n$ coprime to $k$.
\end{abstract}


\section{Introduction}

In 1981 Spiro \cite{spiro1981thesis} show that $d(n)=d(n+B)$ with $B=5040$, where $d(n)$ denotes the number of divisor of $n$. Subsequently Heath-Brown \cite{heath1984divisor} modify spiro's method to solve the case $B=1$. In 1997, Pinner \cite{pinner1997repeated} adapt method of Heath-Brown to deal with an arbitrary $B$.

In this paper, we show similar results for some restricted divisor functions. Define

\begin{equation}
    \notag
    d_k(n)=\sum_{\emptyfrac{m|n}{(m,k)=1}}1.
\end{equation}

\noindent
Evidently $d_k(n)$ is multiplicative. We have

\begin{theorem*}
For any positive integer $k$ and $B$, $d_k(n)=d_k(n+B)$ has infinitely many solutions. Moreover, the number of solutions $n\leq x$ is $\gg x/(\log x)^7$.
\end{theorem*}

As in Heath-Brown and Pinner, we need the following lemma to construct solutions of $d_k(n)=d_k(n+B)$.

\begin{lemma}
\label{key lemma}
For any positive integers $k,B,N$, there exist $N$ distinct positive integers $a_0,\cdots,a_{N-1}$ with the following properties. Let $d_{mn}=a_m-a_n$, then
\begin{equation}
    \label{gcd}
    d_{mn}\,|\,(a_m,a_n),\ m\neq n,
\end{equation}

\noindent
and
\begin{equation}
    \label{second condition of key lemma}
    d_k(Ba_n)d_k\bracket{\frac{Ba_m}{d_{mn}}}=d_k(Ba_m)d_k\bracket{\frac{Ba_n}{d_{mn}}},\ m\neq n.
\end{equation}
\end{lemma}

~

\section{Proof of Lemma \ref{key lemma}}

The proof follows Heath-Brown and Pinner. In their original proofs, the certain power of the first $N$ primes are arranged. However, these primes may not contribute to $d_k(n)$. We cannot remove the primes power from $a_n$ since they assure the condition \eqref{gcd}, as discussed in \cite{heath1984divisor}. We deal with this by appending other primes power in $a_n$. It is the new ingredient in this paper.

Without loss of generality, we assume $k$ is squarefree and $k\geq2$. It suffices to prove Lemma \ref{key lemma} for $N=2^s\geq k$. Let $\omega=\omega(k)$ and $G=\bbZ_2^s$ and set
\begin{equation}
    \notag
    n(\pi)=\sum_{j=0}^{s-1}\pi_j2^j,\ \pi=(\pi_0,\cdots,\pi_{s-1})\in G,
\end{equation}

\noindent
then $n$ gives a 1-1 correspondence between $G$ and $\{0,\cdots,N-1\}$. Let $p_\pi=p_{n(\pi)+1},\ \pi\in G$ denote the first $N$ primes. Evidently $k\,|\,\prod p_\pi$. Let
\begin{equation}
    \notag
    k=p_{\mu_1}\cdots p_{\mu_\omega},\ \mu_1,\cdots,\mu_\omega\in G,\ p_{\mu_1}<\cdots<p_{\mu_\omega},
\end{equation}

\noindent
and set
\begin{equation}
    \notag
    G^*=G\backslash\{\mu_1,\cdots,\mu_\omega\},\ P=\prod_{1\leq i\leq N+\omega}p_i.
\end{equation}

\noindent
Assuming that
\begin{equation}
    \notag
    B=B_1\prod_{\pi\in G}p_\pi^{\alpha_\pi}\prod_{1\leq i\leq\omega}p_{N+i}^{\alpha_{N+i}},\ (B_1,P)=1.
\end{equation}

\noindent
We shall arrange that
\begin{equation}
    \label{power of first N primes}
    p_\pi^{(\alpha_\pi+1)n(\pi+\sigma)}\,||\,a_{n(\sigma)},\ \pi,\sigma\in G,
\end{equation}
\begin{equation}
    \label{power of N+1 to N+w primes}
    p_{N+i}^{(\alpha_{N+i}+1)n(\mu_i+\sigma)}\,||\,a_{n(\sigma)},\ 1\leq i\leq\omega,\ \sigma\in G,
\end{equation}

\noindent
and for the remaining primes $p\nmid P$ we will prove
\begin{equation}
    \label{|->||}
    p\,|\,d_{n(\sigma)n(\tau)}\Rightarrow p\,||\,d_{n(\sigma)n(\tau)}\text{ and }p\,||\,a_{n(\sigma)}.
\end{equation}

\noindent
Since $n(\pi+\sigma)\neq n(\pi+\tau)$ whenever $\sigma\neq\tau$, and together with \eqref{|->||} giving \eqref{gcd}. The first $N+\omega$ primes together contribute
\begin{equation}
    \label{contribution of first N+w primes}
    \begin{aligned}
        &\ \,\,\prod_{\pi\in G^*}(\alpha_\pi+1)^2(n(\pi+\sigma)+1)(n(\pi+\tau)+1-\min\{n(\pi+\sigma),n(\pi+\tau)\})\\
        &\cdot\prod_{1\leq i\leq\omega}(\alpha_{N+i}+1)^2(n(\mu_i+\sigma)+1)(n(\mu_i+\tau)+1-\min\{n(\mu_i+\sigma),n(\mu_i+\tau)\})
    \end{aligned}
\end{equation}

\noindent
to $d_k(Ba_{n(\sigma)})d_k(Ba_{n(\tau)}/d_{n(\sigma)n(\tau)})$, since those primes $p_{\mu_i}$ do not contribute. Then we replace $\pi+\sigma$ by $\rho$ in \eqref{contribution of first N+w primes} to obtain
\begin{equation}
    \notag
    \prod_{1\leq i\leq\omega}\frac{(\alpha_{N+i}+1)^2}{(\alpha_{\mu_i}+1)^2}\prod_{\rho\in G}(\alpha_\rho+1)^2(n(\rho)+1)(n(\rho+\sigma+\tau)+1-\min\{n(\rho),n(\rho+\sigma+\tau)\}).
\end{equation}

\noindent
Since this is symmetrical in $\sigma$ and $\tau$, the corresponding factors on both sides of \eqref{second condition of key lemma} agree. For the remaining primes $p\nmid P$ we have $(p,k)=1$. Suppose that $p^e\,||\,B_1,\ e\geq0$. If $p\,|\,d_{n(\sigma)n(\tau)}$, then by \eqref{|->||} we obtain
\begin{equation}
    \notag
    p^{e+1}\,||\,Ba_{n(\sigma)},\ p^e\,||\,Ba_{n(\tau)}/d_{n(\sigma)n(\tau)},\ p^{e+1}\,||\,Ba_{n(\tau)},\ p^e\,||\,Ba_{n(\sigma)}/d_{n(\sigma)n(\tau)};
\end{equation}

\noindent
if $p\nmid d_{n(\sigma)n(\tau)}$, without loss of generality we may assume $p\nmid a_{n(\sigma)}$ and $p^f\,||\,a_{n(\tau)},\ f\geq0$, so
\begin{equation}
    \notag
    p^e\,||\,Ba_{n(\sigma)},\ p^{e+f}\,||\,Ba_{n(\tau)}/d_{n(\sigma)n(\tau)},\ p^{e+f}\,||\,Ba_{n(\tau)},\ p^e\,||\,Ba_{n(\sigma)}/d_{n(\sigma)n(\tau)}.
\end{equation}

\noindent
Thus \eqref{second condition of key lemma} holds in all cases. Now we remain to show how \eqref{power of first N primes}, \eqref{power of N+1 to N+w primes} and \eqref{|->||} can be achieved.

Let $I=(0,\cdots,0)\in G$, we will construct distinct integers $\delta_\sigma,\ \sigma\in G$, such that $\delta_I=0$ and
\begin{equation}
    \notag
    \delta_\sigma-\delta_\tau=E_{\sigma\tau}^*F_{\sigma\tau}^*,\ \sigma,\tau\in G,
\end{equation}

\noindent
and 
\begin{equation}
    \label{E*}
    E_{\sigma\tau}^*=\prod_{\pi\in G}p_\pi^{(\alpha_\pi+1)\min\{n(\pi+\sigma),n(\pi+\tau)\}}\prod_{1\leq i\leq\omega}p_{N+i}^{(\alpha_{\mu_i}+1)\min\{n(\mu_i+\sigma),n(\mu_i+\tau)\}},\ \sigma\neq\tau,
\end{equation}
\begin{equation}
    \label{(F*,P)=1}
    (F_{\sigma\tau}^*,P)=1,
\end{equation}
\begin{equation}
    \label{F* is squarefree}
    F_{\sigma\tau}^*\text{ is squarefree},
\end{equation}
\begin{equation}
    \label{F* are coprime}
    (F_{\sigma\tau}^*,F_{\pi\rho}^*)=1\text{ for }\{\sigma,\tau\}\neq\{\pi,\rho\}.
\end{equation}

\noindent
Writing $J=(1,\cdots,1)\in G$, we solve the simultaneous congruences
\begin{equation}
    \label{mod p_pi}
    \modulo{x}{-\delta_{\pi+J}+p_\pi^{(\alpha_\pi+1)(N-1)}}{p_\pi^{(\alpha_\pi+1)(N-1)+1}},\ \pi\in G,
\end{equation}
\begin{equation}
    \label{mod p_(N+i)}
    \modulo{x}{-\delta_{\mu_i+J}+p_{N+i}^{(\alpha_{\mu_i}+1)(N-1)}}{p_{N+i}^{(\alpha_{\mu_i}+1)(N-1)+1}},\ 1\leq i\leq\omega,
\end{equation}
\begin{equation}
    \label{mod F*^2}
    \modulo{x}{-\delta_\sigma-E_{\sigma\tau}^*F_{\sigma\tau}^*}{F_{\sigma\tau}^{*2}},\ \sigma,\tau\in G,\ n(\sigma)<n(\tau).
\end{equation}

\noindent
We take a solution $x$ large enough such that $a_{n(\sigma)}=x+\delta_\sigma$ is positive for all $\sigma\in G$.

From \eqref{mod p_pi} we have
\begin{equation}
    \notag
    p_\pi^{(\alpha_\pi+1)(N-1)}\,||\,x+\delta_{\pi+J}\Rightarrow p_\pi^{(\alpha_\pi+1)n(J)}\,||\,a_{n(\pi+J)};
\end{equation}

\noindent
if $\sigma\neq\pi+J$, then $n(\pi+\sigma)<N-1=n(\pi+(\pi+J))$. From \eqref{E*} and \eqref{mod p_pi} we have
\begin{equation}
    \notag
    p_\pi^{(\alpha_\pi+1)n(\pi+\sigma)}\,||\,a_{n(\sigma)}=(x+\delta_{\pi+J})+E_{\sigma,\pi+J}^*F_{\sigma,\pi+J}^*.
\end{equation}

\noindent
From \eqref{mod p_(N+i)} we have
\begin{equation}
    \notag
    p_{N+i}^{(\alpha_{\mu_i}+1)(N-1)}\,||\,x+\delta_{{\mu_i}+J}\Rightarrow p_{N+i}^{(\alpha_{\mu_i}+1)n(J)}\,||\,a_{n(\mu_i+J)};
\end{equation}

\noindent
if $\sigma\neq\mu_i+J$, then $n(\mu_i+\sigma)<N-1=n(\mu_i+(\mu_i+J))$. From \eqref{E*} and \eqref{mod p_(N+i)} we have
\begin{equation}
    \notag
    p_{N+i}^{(\alpha_{\mu_i}+1)n(\mu_i+\sigma)}\,||\,a_{n(\sigma)}=(x+\delta_{\mu_i+J})+E_{\sigma,\mu_i+J}^*F_{\sigma,\mu_i+J}^*.
\end{equation}

\noindent
Now we verify all cases of conditions \eqref{power of first N primes} and \eqref{power of N+1 to N+w primes}. If $p\nmid P$ and $p\,|\,d_{n(\sigma)n(\tau)}=|E_{\sigma\tau}^*F_{\sigma\tau}^*|$ then $p\,||\,d_{n(\sigma)n(\tau)}$ from \eqref{F* is squarefree} and $p\,||\,a_{n(\sigma)}$ from \eqref{mod F*^2} provided $n(\sigma)<n(\tau)$. If $n(\sigma)>n(\tau)$, then we have
\begin{equation}
    \notag
    a_{n(\sigma)}=a_{n(\tau)}-E_{\tau\sigma}^*F_{\tau\sigma}^*=\modulo{x+\delta_\tau-E_{\tau\sigma}^*F_{\tau\sigma}^*}{-2E_{\tau\sigma}^*F_{\tau\sigma}^*}{F_{\tau\sigma}^{*2}}.
\end{equation}

\noindent
We have $p\,||\,a_{n(\sigma)}$ in this case too, since $p>2$. Now we verify the condition \eqref{|->||}.

Next we construct the distinct integers $\delta_\sigma$ which satisfy \eqref{E*}-\eqref{F* are coprime}. We let $\alpha_{n(\pi)}=\alpha_\pi$ and
\begin{equation}
    \notag
    \alpha=\max_{1\leq i\leq N+\omega}\ \alpha_i.
\end{equation}

\noindent
Define
\begin{equation}
    \notag
    c_\sigma=\prod_{\pi\in G}p_\pi^{(\alpha_\pi+1)n(\pi+\sigma)}\prod_{1\leq i\leq\omega}p_{N+i}^{(\alpha_{\mu_i}+1)n(\mu_i+\sigma)}.
\end{equation}

\noindent
Then we set $\delta_\sigma=c_\sigma-c_I+\beta_\sigma P^{(\alpha+1)N}$ and $\beta_I=0$. Thus $\delta_I=0$. Moreover, since $(\alpha_\pi+1)n(\sigma+\pi)<(\alpha+1)N$, the integers $\delta_\sigma+c_I$ are distinct, so are $\delta_\sigma$. Additionally, \eqref{E*} and \eqref{(F*,P)=1} hold.

It remains to show that there are integers $\beta_\sigma$ such that \eqref{F* is squarefree} and \eqref{F* are coprime} hold. It suffices to show that $p^2\nmid\prod_{n(\sigma)<n(\tau)}(\delta_\sigma-\delta_\tau)$ for all primes $p\nmid P$. We write $c_{n(\pi)}=c_\pi,\ \beta_{n(\pi)}=\beta_\pi$ and define
\begin{equation}
    \notag
    g_M=g_M(\beta_0,\cdots,\beta_{M-1})=\prod_{0\leq n<m\leq M-1}(c_n-c_m+(\beta_n-\beta_m)P^{(\alpha+1)N}).
\end{equation}

\noindent
Since $g_N=\prod_{n(\sigma)<n(\tau)}(\delta_\sigma-\delta_\tau)$, it suffices to prove that $p^2\nmid g_M$ for all primes $p\nmid P$ and $M=1,\cdots,N$. We proof by induction on $M$. For $M=1$, $g_1(\beta_0)=1$. For $M>1$, we write $g_M(\beta_0,\cdots,\beta_{M-1})=\gamma h(\beta_{M-1})$, where
\begin{equation}
    \notag
    h(\beta_{M-1})=\prod_{n=0}^{M-2}(c_n-c_{M-1}+(\beta_n-\beta_{M-1})P^{(\alpha+1)N})
\end{equation}

\noindent
and $\gamma=g_{M-1}$. By induction hypothesis, $p^2\nmid\gamma$ for all primes $p\nmid P$. The induction step requires a $\beta^*$ such that $p^2\nmid h(\beta^*)$ for all $p\nmid\gamma P$. Since $p\nmid P$, the congruence $\modulo{h(x)}{0}{p}$ has at most $M-1$ solutions mod $p$. Now since $p>N>M-1$, there exists an $x_p$ (mod $p$) such that $\notmodulo{h(x_p)}{0}{p}$. By Chinese Remainder Theorem we solve the simultaneous congruences
\begin{equation}
    \notag
    \modulo{\beta^*}{x_p}{p}
\end{equation}

\noindent
for all $p\,|\,\gamma$, $p\nmid P$. We take a fixed solution $\beta^*$, then for all $p\,|\,\gamma$, $p\nmid P$, $p\nmid h(\beta^*+\gamma i),\ i\in\bbZ$. If $p\nmid\gamma P$ and $p^2\,|\,h(\beta^*+\gamma i)$ then
\begin{equation}
    \notag
    p^2\,|\,l_n=(c_n-c_{M-1}+(\beta_n-(\beta^*+\gamma i))P^{(\alpha+1)N})
\end{equation}

\noindent
for some $0\leq n\leq M-2$. Notice that $(l_n-l_m)\,|\,\gamma$, hence we cannot have both $p\,|\,l_m$ and $p\,|\,l_n$ for $m\neq n$. Since for each $n$, there is exactly one $i$ mod $p^2$ such that $p\,|\,l_n$, thus there are at most $M-1$ such $i$ for which $p^2\,|\,h(\beta^*+\gamma i)$. On the other hand, if $p^2\,|\,h(\beta^*+\gamma i)$ for some $0<i\leq x$, then $p^2\leq cx$ (the constant $c$ is independent of $x$). Therefore for sufficiently large $x$, we have

\begin{equation}
    \notag
    \begin{aligned}
    \sum_{p\nmid\gamma P}\#\{i:0<i\leq x,\ p^2\,|\,h(\beta^*+\gamma i)\}&\leq\sum_{N<p\leq(cx)^{1/2}}\ceil{\frac{x}{p^2}}(M-1)\\
    &<N\sum_{N<p\leq(cx)^{1/2}}\bracket{\frac{x}{p^2}+1}\\
    &\leq Nx\sum_{p>N}\frac{1}{p^2}+N(cx)^{1/2}\\
    &<Nx\sum_{i=N+1}^\infty\frac{1}{i^2}\\
    &<x.
    \end{aligned}
\end{equation}

\noindent
Thus for sufficiently large $x$, we can find an $i$ such that $0<i\leq x$ and $p^2\nmid h(\beta^*+\gamma i)$ for all $p\nmid\gamma P$. Hence the integers $\beta_\sigma$ are found and the proof of Lemma \ref{key lemma} is complete.

\section{Proof of the Theorem}

Let $a_n$, $d_{mn}$ be as in Lemma \eqref{key lemma} and define
\begin{equation}
    \notag
    A=BkN!\prod_{0\leq n\leq N-1}a_n.
\end{equation}

\noindent
We choose $N$ distinct primes $p_1,\cdots,p_N$, none of which divides $A$ and we write
\begin{equation}
    \notag
    r_n=p_n^{d_k(Ba_n)-1},
\end{equation}

\noindent
so we have $d_k(r_n)=d(r_n)=d_k(Ba_n)$. By the Chinese Remainder Theorem we solve the simultaneous congruences
\begin{equation}
    \notag
    \modulo{a_nAx+1}{r_n}{r_n^2}.
\end{equation}

\noindent
Let $X$ be a fixed solution and write
\begin{equation}
    \label{definition of Y_n}
    a_nAX+1=r_nY_n,
\end{equation}

\noindent
so
\begin{equation}
    \label{(Y_n,Ar_n)=1}
    (Y_n,Ar_n)=1.
\end{equation}

\noindent
Now we define $R=\Pi r_n$, $R_n=R/r_n$ and
\begin{equation}
    \label{definition of F_n(x)}
    F_n(x)=a_nAR_nRx+Y_n.
\end{equation}

\noindent
We first prove that
\begin{equation}
    \label{(F_n(x),AR)=1}
    (F_n(x),AR)=1.
\end{equation}

\noindent
If $p\,|\,Ar_n$ then \eqref{(Y_n,Ar_n)=1} and \eqref{definition of F_n(x)} imply that $p\nmid F_n(x)$. If $m\neq n$ and $p_m\,|\,F_n(x)$, then $p_m\,|\,Y_n$ by \eqref{definition of F_n(x)}. Hence $(p_m,AX)=1$ from \eqref{definition of Y_n} and again by \eqref{definition of Y_n} we obtain
\begin{equation}
    a_nAX+1\equiv\modulo{a_mAX+1}{0}{p_m}.
\end{equation}

\noindent
It follows that $p_m\,|\,d_{mn}$ so $p_m\,|\,a_m$, which contradicts $(p_m,AX)=1$. This complete the proof of \eqref{(F_n(x),AR)=1}.

From \eqref{definition of Y_n} and \eqref{definition of F_n(x)} we have, if $a_m>a_n$,
\begin{equation}
    \frac{Ba_m}{d_{mn}}r_nF_n(x)=B+\frac{Ba_n}{d_{mn}}r_mF_m(x).
\end{equation}

\noindent
Moreover, $Ba_m/d_{mn}$, $r_n$ and $F_n(x)$ are coprime in pairs, by \eqref{(F_n(x),AR)=1}. Hence
\begin{equation}
    d_k\bracket{\frac{Ba_m}{d_{mn}}r_nF_n(x)}=d_k\bracket{\frac{Ba_m}{d_{mn}}}d_k(Ba_n)d_k(F_n(x)).
\end{equation}

\noindent
Similarly
\begin{equation}
    d_k\bracket{\frac{Ba_n}{d_{mn}}r_mF_m(x)}=d_k\bracket{\frac{Ba_n}{d_{mn}}}d_k(Ba_m)d_k(F_m(x)).
\end{equation}

\noindent
By Lemma \ref{key lemma} we have $d_k(Ba_n)d_k(Ba_m/d_{mn})=d_k(Ba_m)d_k(Ba_n/d_{mn})$, so $d_k(u)=d_k(u+B)$ with
\begin{equation}
    \notag
    u=\frac{Ba_n}{d_{mn}}r_mF_m(x),
\end{equation}

\noindent
provided
\begin{equation}
    \label{d_k(F_m(x))=d_k(F_n(x))}
    d_k(F_m(x))=d_k(F_n(x)).
\end{equation}

\noindent
By \eqref{(F_n(x),AR)=1} we obtain $d_k(F_n(x))=d(F_n(x))$, thus \eqref{d_k(F_m(x))=d_k(F_n(x))} is equivalent to
\begin{equation}
    d(F_m(x))=d(F_n(x)).
\end{equation}

\noindent
If $a_m<a_n$, one can merely replace $m$ and $n$. The rest of the proof is closely follows Heath-Brown \cite{heath1984divisor}. For convenience we sketch it here. We need the following lemma (Theorem 10.5 of \cite{halberstam2013sieve}).
\begin{lemma}
\label{N-dimensional sieve}
Let $g$ be a natural number $>1$, and let $a_i,b_i(i=1,\cdots,g)$ be integers satisfying
\begin{equation}
    \prod_{i=1}^ga_i\prod_{1\leq t<s\leq g}(a_tb_s-a_sb_t)\neq0.
\end{equation}

\noindent
Suppose also that
\begin{equation}
    \prod_{i=1}^g(a_in+b_i)
\end{equation}

\noindent
has no fixed prime divisor. Then, for any natural number $r$ satisfying
\begin{equation}
    r>(g+1)\log v_g+g-1,
\end{equation}

\noindent
there is a positive number $\delta$ such that, as $x\rightarrow\infty$,
\begin{equation}
    \label{lower bound}
    \#\largebracket{n:1\leq n\leq x,\prod_{i=1}^g(a_in+b_i)=P_r}\geq\delta\frac{x}{(\log x)^g}(1+O((\log x)^{-1/2})),
\end{equation}

\noindent
where $\delta$ and the $O$ constant depend only on $r,g$ and on the $a_i$'s and $b_i$'s.
\end{lemma}

\noindent
As is shown by Xie \cite{xie1983k}, we can take the pair $(g,r)=(7,27)$. Moreover, according to \cite{heath1984divisor}, we can require the numbers $a_in+b_i$ to be squarefree, without affecting the lower bound in \eqref{lower bound}. Then we shall verify all $F_n(x)$ satisfy conditions of Lemma \ref{N-dimensional sieve} and conclude that there are $\gg x/(\log x)^7$ solutions of $d(F_m(x))=d(F_n(x))$.


\end{document}